\documentclass{ifacconf}

\usepackage{natbib}            % you should have natbib.sty
\usepackage{graphicx}          % Include this line if your 
\newcommand{\eq}{\begin{equation}\begin{array}{rllllllllllllllllllllllllllllllll}}
\newcommand{\ee}{\end{array}\end{equation}}
\newcommand{\bmt}{\left[ \begin{array}{ccccccccc}}
\newcommand{\emt}{\end{array}\right]}

\begin{document}

\begin{frontmatter}

\title{Model Predictive Regulation\thanksref{footnoteinfo}} % Title, preferably not more than 10 words.

\thanks[footnoteinfo]{This work was supported in part by the AFOSR}

\author[First]{Cesar O. Aguilar} 
\author[Second]{Arthur J. Krener} 
%\author[Third]{Third C. Author}

\address[First]{California State University, Bakersfield, CA, 93311, USA\\ (email: caguilar24@csub.edu).}                                              
\address[Second]{University of California, Davis, CA 95616, USA\\ (e-mail: ajkrener@ucdavis.edu)}
          

\begin{abstract}                          % Abstract of not more than 250 words.
We show how optimal nonlinear regulation can be achieved in a model predictive control fashion.
\end{abstract}

\end{frontmatter}

\section{Introduction}
Nonlinear regulation is the design of a feedforward and feedback control law that regulates the output of a nonlinear plant affected by a nonlinear exosystem.  The standard approach involves 
two  steps.  First calculate a feedforward control law that achieves regulation when the states of the plant and exosystem are on a given manifold, called the tracking manifold.  In discrete time this requires the off-line solution of the Francis-Byrnes-Isidori  (FBI) equations.   The second step is to design an additional feedforward and feedback control law that  drives the combined system to the tracking manifold.  This second step can be accomplished by soving off-line the dynamic programming equations (DP) for  an infinite horizon optimal control problem.
Solving either the FBI or the DP equations is difficult when state dimensions are greater than one.
This paper introduces the Model Prediction Regulation (MPR)  method for  nonlinear  regulation.  We show that by choosing a suitable running cost, regulation can be achieved by Model Predictive Control (MPC) type methods, i.e., solving on-line finite horizon nonlinear programs.  We also show how to approximate the ideal terminal cost for the nonlinear program.
MPC is a way of stabilizing a plant to an operating point.  MPR is a way of stabilizing a plant, to an operating trajectory. 

\section{Stabilization around an Operating Point}
\setcounter{equation}{0}

Consider a controlled dynamical system in discrete time 
\eq \nonumber
x^+(t)&=& x(t+1)= f(x(t),u(t))
\ee
where $ x,\ u $ are $ n,\ m$ dimensional.   An equilibrium or operating point
is a pair $(x^e,u^e)$ where 
$
x^e= f(x^e,u^e)
$.
A typical engineering goal is to find a feedback law $u=\kappa(x)$ which makes
the equilibrium state asymptotically stable under the closed loop dynamics,
\eq \nonumber
x^+&=& f(x,\kappa(x))
\ee
at least locally in some neighborhood of $x^e$.
There may also be state and/or control constraints of the form
\eq \nonumber
0&\le& g(x,u)
\ee
which compound the problem.  We shall assume that these constraints are not active in a neighborhood of $(x^e,u^e)$.
By translations of the state and the control coordinates
one can assume 
 that $ x^e=0,\ u^e=0$.

Stabilization is a difficult problem in part because it does not have an unique solution.  
It is generally easier to solve a problem that has an unique solution.
Therefore a standard approach is to recast the stabilization problem as an infinite horizon optimal control problem,
\eq \nonumber
\min_u \sum_{t=0}^\infty\ l(x(t),u(t)) 
\ee
subject to 
\eq \nonumber
x^+= f(x,u), &\quad\ &
0\le g(x,u)
\ee

The Lagrangian or running cost $l(x,u)$ may be given by economic considerations or chosen 
so that $ x(t)\to 0$ without using too much $u(t)$ and without violating the constraints, e.g.
 $
l(x,u) =x'Qx+u'Ru
$
where $Q\ge 0$ and $R>0$. It is sometimes useful to incorporate the constraints into   the running cost  by redefining
\eq \nonumber
l(x,u) &=& \left\{ \begin{array}{ccc} x'Qx+u'Ru &\mbox{ if }&  g(x,u)\ge 0\\
\infty & \mbox{ if }&  g(x,u)< 0 \end{array}\right.
\ee

Let $  \pi(x^0) $, denote the optimal cost and $u=  \kappa(x^0) $ denote the optimal control  given $ x(0)=x^0$.
Then it is well-known that they are the solution to the  dynamic programming (DP) equations
\eq \nonumber
\pi(x)&=& \mbox{min}_u \left\{ \pi(f(x,u)) +l(x,u)\right\}\\
\kappa(x) &=& \mbox{argmin}_u \left\{ \pi(f(x,u)) +l(x,u)\right\}
\ee
But this does not greatly simplify the stabilization problem as the DP equations are notoriously difficult to solve if the state dimension $ n $ is greater than $ 2$.  

But if we can solve the DP equations then
the optimal cost $ \pi(x)$ is a Lyapunov function which verifies the stability of the closed loop system
\eq \nonumber
\pi(x^+(t)) &\le& \pi(x(t))
\ee
provided
\eq \nonumber
\pi(x)  >  0,\quad l(x,u)> 0 &\quad\  & \mbox{ if } x\ne 0
\ee
and some other conditions are satisfied.

An increasingly popular way of solving the stabilization problem is  Model Predictive Control (MPC).
Instead of solving the infinite horizon optimal control problem off-line for every possible initial state $x^0$, we solve a finite horizon optimal control problem on-line for the current  state of the system.
If $ x(t)=x^t $ then we pose the finite horizon  optimal control problem
\eq \nonumber
\min_u \sum_{s=t}^{t+T-1} l(x(s),u(s)) +\pi^T(x(t+T))
\ee
subject to 
\eq \nonumber
x^+= f(x,u),\quad\  & 
0\le g(x,u),\quad\  &
x(t)=x^t
\ee
The terminal cost $ \pi^T(x)$ may only be defined in some compact neighborhood ${\cal X} $ of the equilibrium state, $x=0$, so an extra constraint is needed,
\eq \nonumber
x(t+T) &\in &{\cal X}
\ee
We shall assume that the constraints  are not active on ${\cal X}$.

This is a nonlinear program and a fast solver is used  to obtain the optimal control sequence
$
u^*(t), \ldots, u^*(t+T-1)
$.
Then the feedback 
$
\kappa(x(t))= u^*(t)
$
is used for one time step.
The process is repeated at subsequent times.

The key issues in MPC are the following.
\begin{itemize}
\item Since the discrete time model is probably an approximation to a continuous time dynamics, the time step must be short compared to the dynamics.
 \item The horizon $ T $ must be short enough and the nonlinear program simple enough to  be solved in one time step.
  \item The horizon $ T $ must be long enough and/or ${\cal X}$ large enough to ensure the constraint  $ x(t+T)\in{\cal X}$ can be met.
  \item The terminal cost must be a control Lyapunov function for the dynamics, i.e.
\eq \nonumber
  \min_u \pi^T(f(x,u))< \pi^T(x)  
  \ee
  for all $x\in {\cal X}$. 
  \item The ideal 
  terminal cost is the optimal cost of the infinite horizon optimal control problem provided that it can be computed on a large enough ${\cal X}$.  Then the exact solutions to the finite horizon and infinite horizon optimal control problems are identical.
 \end{itemize}
 
\section{Regulation } 
In the regulation problem we are given a plant
\eq \label{pl}
x^+ = f(x,u,w),& \quad\ &
y= h(x,u,w)
\ee
that it is affected by an external signal $ w(t) $ that might be a command or a disturbance.
The dimension of $ y $ is $ p $ and we usually assume that the system is square,  $ p=m$, the dimension of the control.
The goal is to find a feedforward and feedback $ u=\kappa(x,w) $ such that $ y(t)\to 0$ as $ t\to \infty$.

 \cite{Fr77} solved the linear problem assuming that the plant is linear and the external signal is generated by a linear exosystem of the form
\eq \label{ex}
w^+&=&Aw\ee
The dimension of $ w $ is $ k $.

 \cite{IB90}  solved the problem for a smooth nonlinear plant assuming that the external signal is generated by a smooth nonlinear exosystem of the form
\eq \label{ex}
w^+&=& a(w)
\ee
See also  \cite{HR90}.

A usual assumption is that the exosytem is neutrally stable in some sense, e.g., all the eigenvalues of 
\eq \label{A}
A&=& \frac{\partial a} {\partial w} (0)
\ee
are on the unit circle.

The first step in  nonlinear regulation is to solve the discrete  time Francis Byrnes Isidori (FBI) equations, see \cite{HL93}. We wish to find functions $ x=\theta(w) $ and $ u=\alpha(w) $ such that
\eq \label{dfbi}
f(\theta(w),\alpha(w),w) &=& \theta(a(w))\\
h(\theta(w),\alpha(w),w) &=& 0
\ee

Then the graph of $ x=\theta(w) $ is an invariant submanifold of $ (x,w) $ space under the  feedforward control law $u=\alpha(w)$ and on this submanifold  $ y=0$. The next step is  to find an additional feedforward and feedback control law that makes this invariant manifold attractive.

For the continuous time regulation problem \cite{Kr92}, [2001]  cast this next  step  as an optimal control problem in the transverse state and control coordinates defined by 
\eq \label{tc}
z= x-\theta(w), &\quad\ &
v= u-\alpha(w)
\ee
He then chose a suitable running cost $ l(z,v,w)$ and showed that  the continuous time  infinte horizon optimal control problem for the combined system  has a local solution despite the fact that it is a nonstandard optimal control problem. There is no control over part $ w$ of the combined state $ (z,w) $ and the dynamics of $w$ is neutrally stable.  In particular, the linear part of the combined system is not stabilizable.

  The optimal cost  $ \rho(z,w) $ of the transverse problem is a Lyapunov function for the $ z$ dynamics under the feedforward and feedback control law $ v=\beta(z, w)$.    
Hence $ y(t) \to 0$ under the combined control law
\eq \nonumber
u&=&\kappa(x,w)=\alpha(w) +\beta(x-\theta(w), w)
\ee
and the optimal cost is
\eq \nonumber
\pi(x,w)&=& \rho(x-\theta(w), w)
\ee

We describe the discrete time analog.  Assume that $\theta(w),\ \alpha(w)$ satisfy the discrete time FBI equations (\ref{dfbi}).  Define the transverse coordinates $z,\ w$ by  (\ref{tc}). Consider the transverse  optimal control problem  
\eq \nonumber
\min_v \sum_{t=0}^\infty\  l(z(t),v(t),w(t)) 
\ee
subject to 
\eq \nonumber
z^+&= &\bar{f}(z,v,w)\\
w^+&= &a(w)\\
y&=&\bar{h}(z,v,w)
\ee
where 
\eq \nonumber 
\bar{f}(z,v,w)&=&f(z+\theta(w),v+\alpha(w),w)-\theta(a(w))\\
\bar{h}(z,v,w)&=&h(z+\theta(w),v+\alpha(w),w)
\ee
Because the discrete time FBI equations (\ref{dfbi})  are satisfied it follows immediately
\eq \label{bar}
\bar{f}(0,0,w)=0,&\quad\ &
\bar{h}(0,0,w)=0
\ee
 
Again this is a nonstandard optimal control problem because there is no control over part $ w$ of the combined state $ (z,w) $. Let $\rho(z,w)$ denote the optimal cost and $\beta(z,w)$ denote
the optimal feedback then they
satisfy  the DP equations 
%%%%%%%
\eq \label{tdp}
\rho(z,w)&=& \mbox{min}_v \left\{ \rho(\bar{f}(z,v,w),a(w)) +l(z,v,w)\right\}\\
\beta(z,w) &=& \mbox{argmin}_v  \left\{ \rho(\bar{f}(z,v,w),a(w)) +l(z,v,w)\right\}
\ee
A solution to these equations  exists locally around $z=0,\ w=0$ if  the linear part of the plant 
\eq \label{plantzv}
\dot{z} &=& Fz+Gv+O(z,v,w)^2\\
y&=& Hz+Jv+O(z,v,w)^2
\ee
and the quadratic part of the running cost
\eq \nonumber
l(z,v,w)&=& z'Qz+2z'Sv+v'Rv+O(x,u,w)^3
\ee
 satisfies the LQR conditions.  The LQR conditions are that $F,G$ is stablizable, $Q^{1/2},F-R^{-1}S'$ is  detectable,
 $Q\ge 0$ and $R>0$. 
 Because of (\ref{bar}) there is no linear terms in $w$  in (\ref{plantzv}) and
\eq \nonumber
 F=\frac{\partial f}{\partial x}(0,0,0),&\quad\ & G=\frac{\partial f}{\partial u}(0,0,0)\\
 \\
 H=\frac{\partial h}{\partial x}(0,0,0),&\quad\ & J=\frac{\partial h}{\partial u}(0,0,0)
   \ee
 
 If $H,F$ is detectable then one choice of the running cost is
\eq \nonumber
l(z,v,w)&=&|\bar{h}(z,v,w)|^2+|v|^2=|y|^2+|v|^2
\ee

Optimal regulation seems to require the off-line solution of both the   FBI  and  the DP equations.  These are two difficult tasks.
Can we use an MPC  approach instead?
The answer is yes provided that the plant is linearly minimum phase which we define below.

\cite{FM13} have proposed a two step MPC approach.  First compute on-line for the current $w(t)$, the ideal plant state and control sequences that  are necessary to keep $y=0$.  Then solve on-line  an optimal control problem that drives the actual plant state and control sequences to the ideal plant and control sequences.  

The method that we are proposing does this in one step.
 The key  is to choose a running cost $ l(x(t),u(t),w(t)) $ that is zero when $ y(t)=0$.  It should also be nonegative definite in $ z(t) $ and positive definite in $ v(t)$ even though we may not know what $ z=x-\theta(w)$, $ v=u-\alpha(w)$ are.
How do we do this?
By making $ l $ a function of $ y(t), y(t+1), \ldots , y(t+r) $ where $ r  $ is the relative degree of the plant.  We elaborate on this  in the following sections.

\section{Relative Degree and Zero Dynamics}
For simplicity of exposition we assume a SISO system, $ m=p=1 $, of the form
\eq \label{sys}
x^+= f(x,u)&,\quad\ & y=h(x,u)
\ee
where $x$ is of dimension $n$.
Define a family of functions $ h^{(j)}(x,u) $ as
\eq \nonumber
h^{(0)}(x,u)= h(x,u), & \quad\ &
h^{(j)}(x,u)= h^{(j-1)}(f(x,u),u)
\ee
The system
has  well-defined relative degree $ r $ if for all $ x,u$
\eq \nonumber
\frac{\partial h^{(j)}}{\partial u}(x,u)\begin{array}{cccccc} &=&0 & \mbox{if} &0\le  j<r\\
&\ne&0 & \mbox{if} &j=r
\end{array}
\ee
In other words $ y(t+r) $ is the first output influenced by $ u(t)$.  

Assuming a  well-defined relative degree $r$ then
$
h^{(j)}(x,u)= h^{(j)}(x)
$
for $0\le j<r$ and 
\eq \nonumber
y(t+j)= h^{(j)}(x(t)),&\quad\ &
y(t+r)= h^{(r)}(x(t),u(t))
\ee
The zero dynamics is obtained by setting $y(t+j)=0$ for $0\le j\le r$.  If there exists a feedback  $u=\gamma(x)$  such that
\eq \label{hr}
0&=&h^{(r)}(x,\gamma( x))
\ee
then the zero dynamics is obtained by closing the loop 
\eq \label{gam}
x^+&=& f(x,\gamma( x))
\ee 
By the well-defined relative degree assumption, $h^{(r)}(0,0)=0$ and $\frac{\partial h^{(r)}}{\partial u}(0,0)\ne0$, so the implicit function theorem applied to (\ref{hr}) implies that $\gamma(x)$ exists in some neighborhood of $x=0$ and $\gamma(0)=0$.

Clearly $u=\gamma( x)$ leaves the set 
\eq \nonumber
\left\{x: h^{(j)}(x)=0,\ 0\le j<r\right\}
\ee
invariant.  It can be shown that at least locally around $x=0$, this set   is a submanifold of $x$ space of dimension $n-r$. 
We call this the zero  manifold of the system (\ref{sys}). 

 The zero dynamics is the closed loop dynamics (\ref{gam}) restricted to the zero manifold.  
The nonlinear  system (\ref{sys}) is said to be minimum phase if the zero dynamics is locally asymptotically stable.  The nonlinear  system (\ref{sys}) is said to be linearly minimum phase if the eigenvalues of the linear part of zero dynamics are strictly inside the unit circle.  The nonlinear system is said to have hyperbolic zero dynamics if the none of the eigenvalues of the linear part of the zero dynamics are on the unit circle.

\section{Infinite Horizon Optimal Regulation } 
We return to the  problem of regulating the plant (\ref{pl}) that is affected by the exosystem (\ref{ex}).
For simplicity of exposition we shall assume a SISO plant, $m=p=1$.

Assume that the poles  of the linear part of exosystem (the eigenvalues of $A$ defined by (\ref{A})) are on the unit circle, the plant with $w=0$ has well-defined relative degree $r$ and the plant has hyperbolic zero dynamics.
Then the discrete time FBI equations (\ref{dfbi}) are solvable,  see \cite{HL93}, \cite{Hu04}. 

As before let 
\eq \nonumber
h^{(0)}(x,u,w)&=& h(x,u,w)\\
h^{(j)}(x,u,w)&=& h^{(j-1)}(f(x,u,w),u,a(w))
\ee
Since $r$ is the relative degree 
\eq \nonumber
h^{(j)}(x,u,w)&=&h^{(j)}(x,w)
\ee
for $0\le j<r$ and there exist locally around $(x,w)=(0,0)$ an unique feedforward and feedback $u=\gamma(x,w)$ such that 
\eq \label{dgam}
0&=&h^{(r)}(x,\gamma( x,w),w)
\ee
The zero manifold of the combined system is
\eq \label{zm}
{\cal Z}&=& \left\{ (x,w): h^{(j)}(x,w) =0,\  0\le j<r\right\}
\ee
and this is invariant under the closed loop dynamics when $u=\gamma(x,w)$.

If $\theta(w),\ \alpha(w) $ satisfy  the discrete FBI equations then for $0\le j<r$ and
\eq \nonumber
h^{(j)}(\theta(w),w)=0, & \quad\ &
h^{(r)}(\theta(w),\alpha(w),w)=0
\ee
so the tracking manifold 
\eq \label{phi}
\left\{ (x,w): x=\theta(w)\right\}
\ee 
is contained in the zero manifold ${\cal Z}$ and $\gamma(x,w)$  is an extension of $\alpha(w) $ off this manifold, 
$
\alpha(w)= \gamma(\theta(w),w)
$

We assume that the zero dynamics of the plant is hyperbolic, none of the eigenvalues of its linear part are on the unit circle.  Because of triangular nature of the combined system, these eigenvalues are inherited by the zero dynamics of the combined system.  The remaining eigenvalues of the linear part of the zero dynamics of the combined system are those of the linear part of the exosystem.  Therefore the manifold (\ref{phi}) is the center manifold of the zero dynamics.

Assuming $r\ge 1$, choose a running cost $l(x,u,w)$ of the form
\eq \label{rc}
l(x,u,w)&=&(h^{(0)}(x,u,w))^2+(h^{(r)}(x,u,w))^2
\ee
and consider the infinite horizon combined optimal control problem 
\eq \nonumber
\min_u \sum_{t=0}^ \infty l(x(t),u(t),w(t)
\ee
subject to
(\ref{pl}) and (\ref{ex}).

Let $\pi(x,w)$ denote the optimal cost and $\kappa(x,w)$ denote the optimal feedback  for this problem.The DP equations are 
\eq \label{csdp}
\pi(x,w)&=& \mbox{min}_u \left\{ \pi(f(x,u,w),a(w)) +l(x,u,w)\right\}\\
\kappa(x,w) &=& \mbox{argmin}_u \left\{ \pi(f(x,u,w),a(w)) +l(x,u,w)\right\}
\ee

It can be shown that if 
\begin{itemize}
\item the poles of the linear part of the  exosystem are on the unit circle,
\item the plant has well-defined relative degree $r$ and is linearly minimum phase,
\item the linear part of the plant is stabilizable
%\item $c_0>0,\ c_r>0$,
\end{itemize}
then
\begin{itemize}
\item the FBI equations (\ref{dfbi}) for regulating (\ref{pl}) and (\ref{ex}) are solvable for $\theta(w),\ \alpha(w)$ locally around $(x,w)=(0,0)$,
%\item the running cost (\ref{rc})  expressed  in transverse coordinates $l(z,v,w) $ is of the form 
%\eq \nonumber
%l(z,v,w)&=& |Hz|^2+|HF^rz+HF^{r-1}Gv|^2\\&&+ O(z,v,w)^3
%\ee
%and
%\eq \nonumber
%0&=&l(0,0,w)\\
%0&=&\frac{\partial l}{\partial (z,v)}(0,0,w)\\ 
%\ee
\item the DP equations (\ref{tdp})  for the transverse optimal control problem are solvable for $\rho(z,w),\ \beta(z,w)$  locally around $(z,w)=(0,0)$,
\item the DP equations (\ref{csdp})  for the combined optimal control problem are solvable and $\pi(x,w)=\rho(x-\theta(w)),\ \kappa(x,w)=\alpha(w)+\beta(x-\theta(w),w)$ locally around $x=\theta(w),w=0$.
\end{itemize}

Notice that to show that  the FBI equations are solvable we assumed that the plant has hyperbolic zero dynamics.   But to show that  the combined optimal contol problem of minimizing (\ref{rc}) subject to (\ref{pl}, \ref{ex}), we rquired the stronger condition that zero dynamics be linearly minimum phase.

If the DP equations for the combined infinite horizon optimal control problem are solvable  for $\pi(x,w),\ \kappa(x,w)$ then the set 
\eq \nonumber
\left\{ (x,w): \pi(x,w)=0\right\}
\ee
is the zero manifold (\ref{zm}) of the combined system.  It is  invariant under the closed loop dynamics $u=\kappa(x,w)$.  Again the tracking manifold (\ref{phi}) is the center manifold of the zero dynamics.

\section{Model Predictive Regulation } 
Standard software cannot solve the DP equations (\ref{csdp})  for the linear-quadratic part of the combined optimal control problem even when the combined dimension is low, e.g.  $n=k=1$. But we can use an MPC approach which we call Model Predictive Regulation (MPR).

Consider a finite horizon version of the combined optimal control problem of the last section.
\eq \nonumber
\min_u \sum_{s=t}^{t+T-1} l(x(s),u(s),w(s)) +\pi^T(x(t+T),w(t+T))
\ee
subject to 
\eq \nonumber
x^+= f(x,u,w)),& \quad\ &
w^+=a(w)\\
y=h(x,u,w), & \quad\ &
0\le  g(x,u,w)\\
x(t)=x^t, &\quad\ &
w(t)=w^t
\ee

We continue to assume that the combined system has a well-defined relative degree $r$ and we choose the running cost as before (\ref{rc}).

The terminal cost $ \pi^T(x,w)$ may only be defined in some compact set ${\cal XW} $ of the combined $(x,w)$ state space so an extra constraint is needed,
\eq \label{XW}
(x(t+T),w(t+T) )&\in &{\cal XW}
\ee
The ideal terminal cost would be the $\pi(x,w) $ of the infinite horizon combined optimal control problem of the previous section.  In the next section we shall show how to approximate this on some compact subset ${\cal XW} $ of the combined state space.  

We shall also assume that the constraints are not active on this compact  set,  $g(x,u,w)>0$ on ${\cal XW}$.

We solve this problem  using a fast solver for the  nonlinear program to obtain the optimal sequence
$
u^*(t), \ldots, u^*(t+T-1)
$.
Then in MPR fashion the  feedforward and  feedback 
\eq \nonumber
\kappa(x(t),w(t))&=& u^*(t)
\ee
is used for one time step.
The process is repeated at subsequent times.

Clearly if the infinite horizon combined optimal control problem has a solution, $ \pi^T(x,w)=\pi(x,w)$ on $ {\cal XW}$, if the constraint (\ref{XW}) can be enforced and if the fast solver delivers the true solution to the finite horizon combined optimal control problem and it is unique then 
  the solutions to the infinite horizon combined optimal control problem and the finite horizon combined optimal control problem are the same.  
  
    It is also desirable to have  the optimal feedback $\kappa(x,w)$ at least on $ {\cal XW}$.  If 
  $
u^*(t), \ldots, u^*(t+T-1)
$
is the optimal control sequence fot the finite horizon optimal control problem at time $t$ and
 $
x^*(t+1), \ldots, x^*(t+T)
$
is the corresponding state trajectory.
Then as initial guess for the problem at time $t+1$ we take  $
u^*(t+1), \ldots, u^*(t+T-1), \kappa(x^*(t+T)
$.   This yields a good initial guess for the next nonlinear program and it speeds up its solution.

  \section{Approximate Solutions}
  
As we have seen the ideal terminal cost $\pi^T(x,w)$ for the finite horizon optimal control problem is the solution $\pi(x,w)$ to the infinite horizon problem.   We would like to approximate the latter on some compact subset ${\cal XW}$ which is forward invariant under the combined dynamics when regulation is being achieved.  If this subset is to be a neighborhood of the origin $(x,w)=(0,0)$, we can proceed as follows.

Assume that Taylor series of  the plant and exosystem around $(x,u,w)=(0,0,0)$ is given by
\eq \nonumber
x^+&=& Fx+Gu+ Bw +f^{[2]}(x,u,w)+\ldots\\
y&=& Hx+Ju+ Dw +h^{[2]}(x,u,w)+\ldots\\
w^+&=& Aw+a^{[2]}(w)+\ldots
\ee
where $f^{[2]}(x,u,w)$ is quadratic vector field in $(x,u,w)$ etc.

The simplest approximation to the optimal cost $\pi(x,w)$ of the infinite horizon optimal control problem  is obtained by first solving the Francis equation for the linear parts of the plant and exosystem.  Find $T,\ L$ such that
\eq \nonumber
\bmt F& G\\H&J\emt \bmt T\\L\emt- \bmt T\\0\emt A&=& -\bmt B\\D\emt
\ee  
  Assuming $m=p$ these equations are solvable if the plant is linearly minimum phase and the poles of $A$ are on the unit circle.

  Then define the approximate transverse coordinates
  \eq \nonumber
  z=x-Tw,&\quad\  & v=u-Lw
  \ee
  After making this change of coordinates, the linear part of the plant becomes
  \eq \label{lplzv}
  z^+= Fz+Gv, & \quad\ &  y= Hz+Jv
  \ee
  Assuming the plant has  relative degree $r>0$ so that $J=0$ then the quadratic part of the running cost (\ref{rc})  is
  \eq \label{lzv}
  l(z,v,w)&=& |Hz|^2+|H\left(F^rz+F^{r-1}Gv\right)|^2\\
  &=& z'Qz+2z'Sv+v'Rv
  \ee
  for suitable choices of $Q,\ R, \ S$.
  
  It is critical that the LQR problem of minimizing the future sum of (\ref{lzv}) subject to (\ref{lplzv}) have a nice solution, that is, a solution where all the closed loop eigenvalues are inside the unit circle.    Define a new contol $\nu=H\left(F^rz+F^{r-1}Gv\right)$ then the LQR problem can be rewritten as 
  \eq \nonumber
  \min_\nu \sum_{t=0}^\infty  |y|^2+|\nu|^2\
  \ee
  subject to 
\eq \label{lplznu}
  z^+=\bar{F} z+\bar{G} \nu, &\quad\   y= Hz
  \ee
  where
  \eq \nonumber
  \bar{F}= F-{HF^r\over HF^{r-1}G}, &\quad\ &
  \bar{G}= {G\over HF^{r-1}G}
   \ee
  This LQR problem has a nice solution if $\bar{F}, \bar{G}$ is stabilzable and $H,\bar{F}$ is detectable.  If $F,G$ is stabilzable then so is $\bar{F}, \bar{G}$ because they differ by state feedback.  If the zero dynamics is linearly minimum phase then $H,\bar{F}$ is detectable.
  Then the optimal cost $z'Pz$ and the optimal feedback $v=Kz$ are given by the familiar LQR equations 
  \eq \nonumber
  P&=& F'PF - K'(R+G'PG)K + Q\\
K&=&-(G'PG+R)^{-1}(G'PF+S')
  \ee
 
 We define the terminal cost   and the terminal control 
 \eq  \label{21}
 \pi^T(x,w)&=&(x-Tw)'P(x-Tw)\\
\kappa^T(x,w)&=&Lw+K(x-Tw)
 \ee
 This is a quadratic-linear approximation to the solution to the infinite horizon optimal control problem. 
 If the solution to the finite horizon optimal control problem at time $t$ is $u^*(t), \ldots, u^*(t+T-1)$  then as an initial guess when computing the solution at time $t+1$, we take $ u^*(t+1), \ldots, u^*(t+T)$ where $u^*(t+T)=\kappa^T(x,w)$.  
 
 We then compute a set $ {\cal XW}$  where $\pi^T(x,w)$ is a Lyapunov function, that is a set of $(x,w)$ where \\
$
 \pi^T(f(x,\kappa^T(x,w),w),a(w))\le \pi^T(x,w)
$.
If this set is large enough so that the terminal constraint (\ref{XW}) can be enforced then we use (\ref{21}) as the terminal cost and feedback of the finite horizon optimal control problem.

Alternatively we can  compute a set $ {\cal XW}$  where $\pi^T(x,w)$ is a control Lyapunov function, that is a set of $(x,w)$ where \\
$
 \min_u \pi^T(f(x,u,w),a(w))\le \pi^T(x,w)
$
and define the terminal feedback by $\kappa^T(x,w)= \mbox{argmin}_u \pi^T(f(x,u,w),a(w))
$.

If  $ {\cal XW}$ is not large enough  then we can increase the horizon $T$ or try to approximate $\pi(x,w)$ more accurately and on a larger set.  

Following \cite{HR92} we can get a higher degree Taylor series approximations to the solution of the discrete FBI equations (\ref{dfbi}).  Due to space limitations we discuss only the extension to degree $2$ approximations but this can be extended to higher degree if the plant and exosystem are smooth enough.  Suppose  
\eq \nonumber
\theta(w)= Tw +\theta^{[2]}(w) ,\quad\ & \alpha(w)= Tw +\alpha^{[2]}(w) 
 \ee
 we plug these expressions into the discrete FBI equations, collect terms of degree $2$ and obtain
 \eq \nonumber
 \bmt F&G\\H&J \emt \bmt \theta^{[2]}(w)\\  \alpha^{[2]}(w) \emt
 -\bmt \theta^2(Aw)\\0\emt \\
 =\bmt Ta ^{[2]}(w)\\0\emt-\bmt f^{[2]}(Tw,Lw,w)\\h^{[2]}(Tw,Lw,w)\emt
 \ee
If $m=p$ these are a square set of linear equations for the unknowns $\theta^{[2]}(w), \ \alpha^{[2]}(w)$.  There is an unique solution if the zero dynamics of the plant is linearly minimum phase and the poles of the exosystem are on the unit circle.

The degree two approximations to the transverse coordinates are
  \eq \nonumber
  z=x-Tw-\theta^{[2]}(w),&\quad\  & v=u-Lw-\alpha^{[2]}(w)
  \ee
  After making this change of coordinates, the linear and quadratic parts of the plant become
  \eq \label{lplzv2}
  z^+&=& Fz+Gv +\bar{f}^{[2]}(z,v,w) \\
  y&=& Hz+Jv +\bar{h}^{[2]}(z,v,w)
  \ee
where  $\bar{f}^{[2]}(0,0,w)=0,\  \bar{h}^{[2]}(0,0,w)=0$.

For this linear quadratic system we choose the running cost similar to  before
\eq \nonumber
l(z(t,)v(t),w(t))=|y(t)|^2+| y(t+r)|^2\\
=z'Qz+2z'Sv+v'Rv
 +l^{[3]}(z,v,w)+\ldots
\ee
for some $l^{[3]}(z,v,w)$ with the property that 
\eq \nonumber
l^{[3]}(0,0,w)=0,  &\quad\ &
\frac{\partial l^{[3]}(}{\partial (z,v)}(0,0,w)=0
%\frac{\partial^2 l^{[3]}(}{\partial (z,v)^2}(0,0,w)&=&0
\ee

Then the cubic approximation to infinite horizon optimal cost and the quadratic approximation to the optimal feedback can be computed using the discrete time version of the method of  \cite{Al61}  as generalized by \cite{Kr92}, [2001] and they can be used as the terminal cost and feedback for the finite horizon optimal control problem.  They are 
\eq \label{32}
\pi^T(x,w)&=& (x-Tw-\theta^{[2]}(w))'P(x-Tw-\theta^{[2]}(w))\\
&& +\pi^{[3]} (x,w)\\
\kappa^T(x,w)&=&=Lw+K(x-Tw-\theta^{[2](w)})+\kappa^{[2]}(x,w)
\ee
where 
$\pi^{[3]} (x,w)= \rho^{[3]}(x-Tw,w)$ and $\kappa^{[2]}(x,w)=\alpha^{[2]}(w) +\beta^{[2]}((x-Tw,w) $ and $\rho^{[3]}(z,w), \  \beta^{[2]}((x-Tw,w) $ are the solution to the linear equations 
\eq \nonumber
\rho^{[3]}(z,w)&=& \rho^{[3]}((F+GK)z,w)+2z'P\bar{f}^{[2]}(z,Kz,w)\\&&+l^{[3]}(z,Kz,,w)\\
0&=& \frac{\partial \pi^{[3]}(}{\partial w}(z,w)G+2z'P\frac{\partial \bar{f}^{[2]}(}{\partial v}(z,Kz,w)\\
&&+
\frac{\partial l^{[3]}(}{\partial v}(z,Kz,w)+2\beta^{[2]}(z,w)R
\ee
These equations are block triangular, the second unknown $\beta^{[2]}(z,w) $ does not appear
in the first equation.  The first equation is uniquely solvable if the eigenvalues of $F+GK$
are strictly inside the unit circle and this will be true if the linear part of the plant is stabilzable and the plant is linearly minimum phase.

Hopefully  the cubic-quadratic approximation (\ref{32}) is a control Lyapunov function on a larger 
set ${\cal XW}$ than the quadatic-linear approximation (\ref{21}).  If not and the plant (\ref{pl})  and exosystem (\ref{ex}) are sufficiently smooth we can go to higher degree approximations  but they  do not always lead to an increase in the size of ${\cal XW}$.  We have found that frequently going to a quartic-cubic approximation yields a larger ${\cal XW}$ than a quadratic-linear approximation.

\section{Examples}
\subsection{Linear Example}
Plant: $ n=3,\ m=1,\ p=1  $ 
 \eq \nonumber
x^+ &=& \bmt 0&1&0\\ 0&0&1\\0&0&0 \emt x+\bmt 0\\1\\0.5\emt u\\
y&=& x_1-w_1
\ee
There are
three plant poles at $ 0 $ and relative degree is
$ r=2$. There is $n-r=1$ plant zero at $-0.5$ and so the plant is minimum phase. 

Exosystem: $ k=2$ 
 \eq \nonumber
w^+&=& \bmt 0&-1 \\ 1&0\emt w
\ee
The two exosystem poles are at $ \pm i$.

 There are no resonances between the plant zero and the exosystem poles   so the Francis equations are solvable. The solution is
  \eq \nonumber
  x_1= w_1,&\quad\ &x_3=-0.2 w_1-0.4w_2\\
  x_2= -w_2,  &\quad\ &
  u=-0.8 w_1+0.4 w_2
  \ee
  
  If we take
 \eq \nonumber
l(x(t),u(t),w(t))&=& (y(t))^2 +(y(t+2))^2\\
&=&  (x_1(t)-w_1(t))^2\\
&& +(x_3(t)+u(t)+w_1(t))^2
\ee 
Then the solution to the DP equations for the infinite horizon optimal control problem 
is
 \eq \nonumber
\pi(x,w)&=& x_1^2-2x_1w_1+x^2_2+2x_2w_2+w_1^2+w_2^2\\
\kappa(x,w)&=&-x_3-w_1
\ee

The zero set $ {\cal Z} $ of $ \pi(x,w) $ is a closed loop invariant three dimensional subspace of $ (x,w)$ space  given by the equations
\eq
0= x_1-w_1,&\quad\ &
0= x_2+w_2
\ee
In the terminology of Wonham and Morse  this is $ {\cal V}^*$, the maximal $A,\ B$ invariant subspace in the kernel of $ C $ for the combined $ x,\ w $  system.  Notice that it is contains the the two dimensional tracking manifold.  The eigenvalues of the closed loop dynamics on $ {\cal Z} $ are the two eigenvalues $\pm i$ of the exosystem and the zero $-0.5$ of the plant.  Hence the closed loop dynamics on $ {\cal Z} $ converges to the tracking manifold.  If the zero dynamics of the plant were unstable then the closed loop dynamics on $ {\cal Z} $ may diverge even if tracking is achieved.  Because the problem is essentially LQR there is no need for MPR techniques.

\subsection{Nonlinear Example}
We start with a continuous time plant that is an asymmetrically damped pendulum
\eq
\dot{x}_1&=&f_1(x,u)=x_2 \\
\dot{x}_2&=&f_2(x,u)=- \sin x_1 -(x_2+x_2^2+x_2^3)+u
\ee
 
 We discretize the unforced dynamics by a third degree  Lie series with time step
 $ t_s=\pi/6, $ so  period of its  linearization without damping is $ 12 $.
 
\eq
x^+&=& F(x)+Gu\\
F(x)&=& x+f(x,0)t_s+L_f(f)(x,0){t_s^2\over 2}+L^2_f(f)(x,0){t_s^3\over 6}\\
G&=& \bmt 0\\ 1\emt
\ee
where
\eq
L_f(g)(x)&=& \frac{\partial g}{\partial x}(x) f(x)
\ee
The 
relative degree is $r=n= 2$ so there are no zeros.

The  output is
\eq
y&=& x_1-w_1
\ee
where $w_1$ is the first coordinate of a discrete time exosystem with period $ 8$ 
\eq
w^+&=& \bmt {\sqrt{2}\over 2}&{-\sqrt{2}\over 2}\\{\sqrt{2}\over 2}&{\sqrt{2}\over 2}
\emt w
\ee
It is a tracking problem.  The periods of the plant and exosystem are different and the plant is asymmetrically damped.
The unforced, undamped plant is slower than the exosystem.

We start the plant and exosystem at
\eq
x(0)=\bmt 0\\0\emt,&&w(0)=\bmt 1\\0\emt
\ee
then $w_1(t) $ is shown in Figure 1.
We set the Lagrangian to be
\eq
l(x,u,w) &=& (h^{(0)}(x,u,w))^2+ (h^{(2)}(x,u,w))^2
\ee 
and compute the power series to degree $3$ 
of the feedfoward and feedback controller.
The tracking error of the cubic and linear
controller are shown in Figure 2.
The cubic error is smaller than the linear error.  The 
steady state average cubic error is $0.0108$
while the steady state average linear error is 
  $0.0499$.  
  
  MPR is a generalization of MPC, it reduces to MPC when the 
  external signal $w(t)=0$.  So using higher degree approximations to
  the optimal cost and optimal feedback as the terminal cost and terminal 
  feedback are useful in MPC.  But they are more useful in MPR because
  the plant is constantly excited by the exosystem and so the higher 
  degree terms have a more pronounced affect.

     \begin{figure}
% \label{fig1}
 \centering
\includegraphics[width=3in]{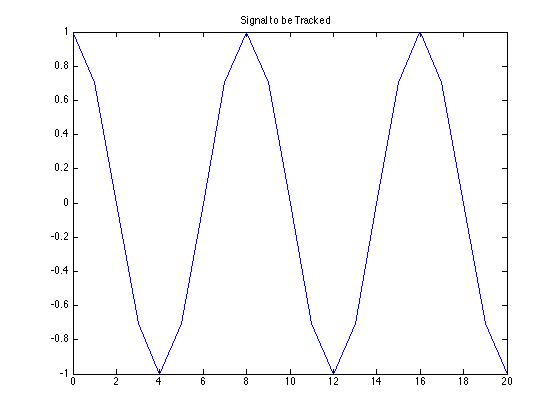}
\caption{Signal to be Tracked}
\end{figure}

     \begin{figure}
% \label{fig1}
 \centering
\includegraphics[width=3in]{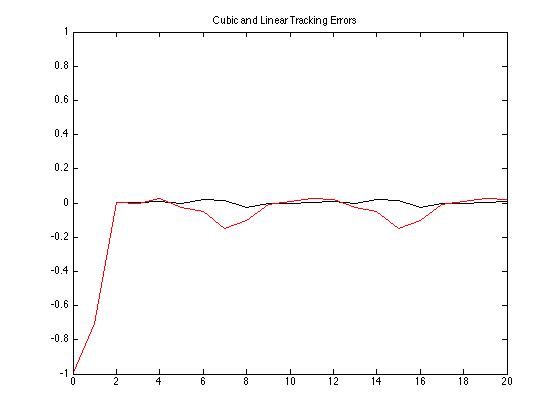}
\caption{Cubic Error in Black, Linear Error in Red}
\end{figure}

If we increase the initial state of the plant to
\eq
x(0)=\bmt 1.5\\0\emt,&&w(0)=\bmt 1\\0\emt
\ee
the cubic controller is able to track but the linear
controller error goes unstable, see Figure 3.

    \begin{figure}
% \label{fig1}
 \centering
\includegraphics[width=3in]{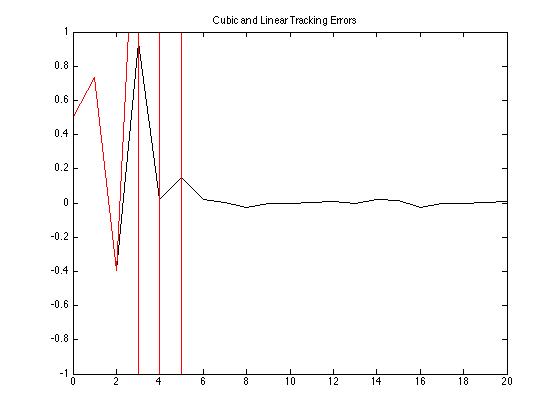}
\caption{Cubic Error in Black, Linear Error in Red}
\end{figure}

Both go unstable when 
\eq
x(0)=\bmt 2\\0\emt,&&w(0)=\bmt 1\\0\emt
\ee
so we need to use MPR.  We set the horizon $T=4$, we use the quartic approximation to $\pi(x,w)$ as the terminal cost and the cubic approximation to $\kappa(x,w)$ as the terminal feedback.  The resulting tracking error is shown in Figure 4 and the control effort is shown in Figure 5.
    \begin{figure}
% \label{fig1}
 \centering
\includegraphics[width=3in]{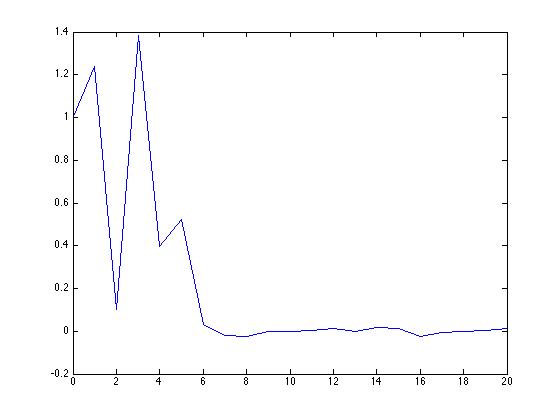}
\caption{MPR Tracking Error}
\end{figure}

    \begin{figure}
% \label{fig1}
 \centering
\includegraphics[width=3in]{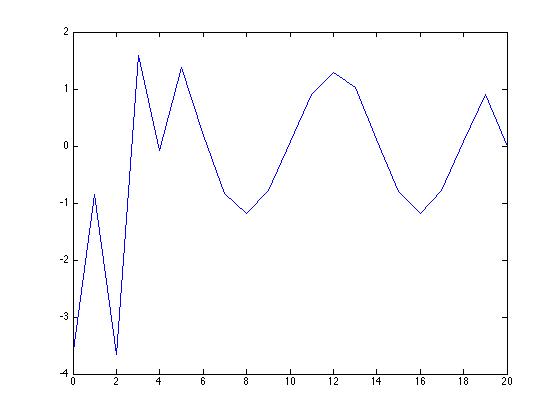}
\caption{MPR Control Effort}
\end{figure}

One of the advantages of MPR is that it can easily handle constraints if they are not active on the tracking manifold. If we impose the constraint $|u(t)|\le 2$ then the resulting tracking error is shown in Figure 6 and the control effort is shown  in Figure 7.  Notice that the  constraint does not affect the tracking performance very much.

  \begin{figure}
% \label{fig1}
 \centering
\includegraphics[width=3in]{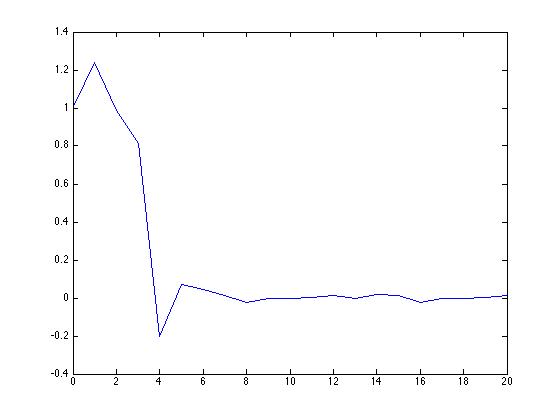}
\caption{Constrained MPR Tracking Error}
\end{figure}

    \begin{figure}
% \label{fig1}
 \centering
\includegraphics[width=3in]{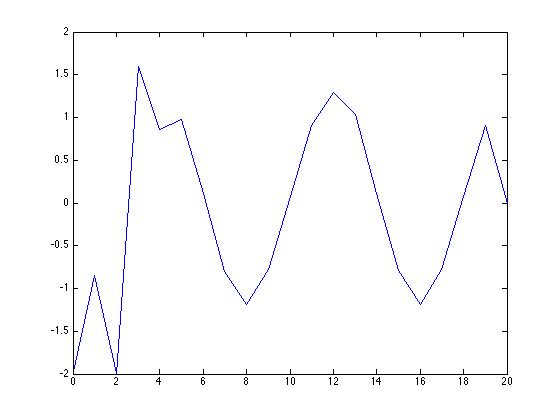}
\caption{Constrained MPR Control Effort}
\end{figure}

Another advantage of MPR is that it is not necessary that the signal $w(t)$ be generated by an exosysytem, all that is needed is that it is known far enough in the future so that the finite horizon optimal control problem makes sense.  In this case one way to generate the terminal cost and terminal feedback is to solve the infinite horizon optimal control problem assuming that the exosystem is  trivial
\eq
w^+&=& w
\ee

\section{Conclusion}  We have shown how nonlinear regulation can be achieved without the off-line solution of the FBI and DP equations by using a MPC approach that we call Model Predictive Regulation (MPR).  MPR requires that the linear part of the plant be  stablizable and minimum phase and that the linear part of the exosystem be neutrally stable.  We have also shown how the infinite horizon optimal cost and feedback can be approximated and used as the terminal cost and feedback in the finite horizon optimal control problem of  MPR.

\begin{ack}                               % Place acknowledgements
The authors are greatful to Professor D. Q. Mayne for suggesting that model predictive control methods could be utilized in nonlinear regulation and for further helpful discussions.
\end{ack}

\end{document}